\documentclass[11pt]{article}
\usepackage{hyperref}
\usepackage{amscd}
\usepackage{amsmath}
\usepackage{amsfonts}
\usepackage{amsthm}
\usepackage{graphicx}
\usepackage{setspace} 
\newtheorem{theorem}{Theorem}

\newtheorem{remark}{Remark}

\begin{document}
\bibliographystyle{plain}
\renewcommand{\a}{\alpha}
\renewcommand{\b}{\beta}
\def\N{\mathbb N}
\def\R{\mathbb R}
\def\Id{{\rm Id}}
\def\mean{\mathop{\bf{E}}}
\def\prob{\mathop{\bf{P}}}
\def\bz{{\bf z}}
\def\rx{{\boldsymbol x}}
\def\bx{X}
\def\rv{{\boldsymbol v}}
\def\bv{V}

\begin{title}
{Robustness of Cucker-Smale flocking model}
\end{title}
\author{
Eduardo Canale\thanks{Instituto de Matem\'atica y Estad\'istica (IMERL), Facultad de Ingenier\'ia, Universidad de la 
Rep\'ublica Julio Herrera y Reissig 565, Montevideo, 11300, Uruguay Email: canale@fing.edu.uy}
\qquad Federico Dalmao\thanks{
Departamento de Matem\'{a}tica y Estad\'{i}stica del Litoral,   
Universidad de la Rep\'{u}blica, Rivera 1350, 50000, Salto, Uruguay. Telephone number: 598 4734 2924. 
E-mail: fdalmao@unorte.edu.uy.} 
\qquad\\ \\ 
Ernesto Mordecki\thanks{
Centro de Matem\'atica, Facultad de Ciencias, Universidad de la Rep\'{u}blica, Igu\'a 4225, 11400, Montevideo, 
Uruguay.
Telephone number: 598 2525 2522,
E-mail: mordecki@cmat.edu.uy.}
\qquad Max O. Souza\thanks{Departamento de Matem\'atica Aplicada, Universidade Federal Fluminense, R. M\'ario 
Santos Braga, s/n, 24020-140, Niter\'oi, RJ.
Telephone number: +55 21 2629 2072
E-mail: msouza@mat.uff.br}
}

\makeatletter
\maketitle
\makeatother
\abstract{
Consider a system of autonomous interacting agents moving in space,
adjusting each own velocity as a weighted mean of the
relative velocities of the other agents. 
In order to test the robustness of the model, we assume that
each pair of agents, at each time step, 
can fail to connect with certain probability, the failure rate.
This is a modification of the (deterministic) Flocking model introduced by Cucker and Smale in \cite{CS06}.
We prove that, if this random failures are independent in time and space, 
and have linear or sub-linear distance dependent rate of decay,
the characteristic
behavior of flocking exhibited  by the original deterministic model, also holds true under random failures,
for all failure rates.
}

\date{}
\makeatletter
\maketitle
\makeatother

\section{Introduction and main result}\label{sec:intro}
The applicability of the results obtained in the mathematical modeling of collective motion,
obviously depend on the accuracy of the assumptions of the model.
Nevertheless, the complexity of the situation of interest requires simplifications,
as very complex models do not admit reasonable mathematical treatment and allow
only simulations or other type of experimental verification. 
In this respect we recommend the
recent review by Vicsek and Zafeiris \cite{v} and its complete list of references.

In this direction, the seminal proposal by Cucker and Smale \cite{CS06}
is an important step in the comprehension of collective motion as, under
reasonable dynamical laws, exhibits cases of flocking, one of the central
issues in collective motion and also non-flocking situations, giving the possibility
of a mathematical treatment of the question.

One of the characteristics of this proposal is that it lets aside some questions,
as the volume of the agents, or its difference in weights, and it also assumes that the
communication between agents is perfect. 
In the present paper we address this third issue.
More precisely, departing from the Cucker-Smale model proposed in \cite{CS06} 
we analyze the situation in which each interaction is subject to random failure. 
Our main result states that the convergence type result obtained in \cite{CS06} holds in our framework, 
with probability one, in the case of agent interactions with linear or sub-linear distance dependent decay's rate, 
obtaining then a robustness result for Cucker-Smale model.
This asymptotic behavior, called \emph{flocking}, consists in the fact that  
for large times, all agent velocities become equal, with fixed 
relative positions. 
A similar approach in the case of \emph{hierarchical} flocking was considered in \cite{dm2} (see also \cite{dm1}).

Consider then a system of $k$ agents with positions and velocities denoted respectively by 	
$\bx=(\bx_1,\dots,\bx_k)$ 
and 
$\bv=(\bv_1,\dots,\bv_k)$. 
All individual positions and velocities are vectors in $\R^3$, and we refer to $\bx$ and $\bv$ as the position
and velocity of the system respectively. 

Assume that the system evolves following the discrete-time dynamic
\begin{equation}\label{eq:cs}
\begin{cases}
\bx_i(t+h)=\bx_i(t)+h\bv_i(t)\\
\bv_i(t+h)=\bv_i(t)+h\sum_{j=1}^ka_{ij}(\bv_j(t)-\bv_i(t)),
\end{cases}
\end{equation}
where $i=1,\dots,k$, $h>0$ is the time step, $a_{ij}(t)$ are the weighting coefficients.

An equivalent way of writing the second equation in \eqref{eq:cs} is
\begin{equation}\label{eq:alternative}  
  \bv_i(t+h)=\left(1-h\sum_{j=1}^ka_{ij}\right)\bv_i(t)+h\sum_{j=1}^ka_{ij}\bv_j(t),
\end{equation}
that, under the natural condition
\begin{equation}\label{eq:h}
0<h\leq {\frac{1}{k}},
\end{equation}
(that makes the first coefficient in the r.h.s. of \eqref{eq:alternative} always positive), 
implies that the velocity at $t+h$ is a linear convex combination of the system velocities at time $t$.
From this follows that the velocity is decreasing in the following sense:
\begin{equation*}
\max_{1\leq i\leq k}\|\bv_i(t+h)\|\leq \max_{1\leq i\leq k}\|\bv_i(t)\|\leq \max_{1\leq i\leq k}\|\bv_i(0)\|,
\end{equation*}
where $\|\bv\|$ is the Euclidean norm of the vector $\bv$ in $\R^3$.
Consequently, in what follows we assume condition \eqref{eq:h}.

In this framework we introduce the possibility that, at each time step, 
each pair of agents $i,j$, can fail to see each other (i.e. they disconnect).
These failures are assumed to be random, independent and with a fixed failure rate probability $\lambda\in (0,1)$.
More precisely, the weighting coefficients in \eqref{eq:cs}, for each pair of agents $j\neq i$ 
are given by 
\begin{equation}\label{eq:coefficients}
a_{ij}(t)=\zeta^{ij}_t\frac{1}{\left(1+\left\|\bx_{i}(t)-\bx_{j}(t)\right\|\right)^{\alpha}},
\end{equation}
where $\zeta^{ij}_t\ (t\geq 0)$ 
are independent and identically distributed Bernoulli random variables with success probability $1-\lambda$, i.e.
$$
\prob(\zeta^{ij}_t=1)=1-\lambda, 
\qquad
\prob(\zeta^{ij}_t=0)=\lambda.
$$
Through the article $\prob$ and $\mean$ denote probability and expectation respectively, and
the (random) event $\zeta^{ij}_t=0$ (resp. $=1$) means that the pair of agents $i,j$ fails to (resp. does) connect
(see each other) at time step $t$.

The factor $\left(1+\left\|\bx_{i}(t)-\bx_{j}(t)\right\|\right)^{-\alpha}$ in \eqref{eq:coefficients}, with $0\leq\alpha\leq 1$,  
is the Cucker-Smale coefficient introduced in \cite{CS06}. 

It is interesting to observe (as was done in \cite{CM07}, and \cite{hl} for continuous time) that the center of mass of the 
system
travels with constant velocity, that happens to be the initial mean velocity of the flock. 
Define
$$
\bar{\bx}(t)=\frac1k\sum_{i=1}^k\bx_i(t),\qquad
\bar{\bv}(t)=\frac1k\sum_{i=1}^k\bv_i(t).
$$
In view of \eqref{eq:cs} (due to the symmetry of the coefficients $a_{ij}=a_{ji}$), 
we obtain that $\bar{\bv}(t+h)=\bar{\bv}(t)$.
We conclude that
\begin{equation}\label{eq:vcxc}
\bar{\bx}(th)=\bar{\bx}(0)+ht\bar{\bv}(0), ,\qquad \bar{\bv}(th)=\bar{\bv}(0),\qquad t=1,2,\dots.
\end{equation}
\begin{theorem}\label{theorem:main}
Consider the dynamical system governed by equations \eqref{eq:cs}, with coefficients
subject to random failure as defined in \eqref{eq:coefficients}. 
\par\noindent {\rm (i)} Under the condition $\alpha<1$ all velocities
tend to a common velocity that is the mean initial velocity $\bar{\bv}(0)$, i.e.
\begin{equation}\label{eq:v-convergence}
(\bv_1(th),\dots,\bv_k(th))\to (\bar{\bv}(0),\dots,\bar{\bv}(0)),\quad\text{as $t\to\infty$.}
\end{equation}
 Furthermore, there exists a limiting configuration $\hat{\rx}=(\hat{\rx}_1,\dots,\hat{\rx}_k)$ that is the limit
 of the relative positions of the system with respect to the center of mass defined in \eqref{eq:vcxc}, 
 more precisely
\begin{equation}\label{eq:x-convergence}
(\bx_1(th)-\bar{\bx}(th),\dots,\bx_k(th)-\bar{\bx}(th))\to (\hat{\rx}_1,\dots,\hat{\rx}_k),\quad\text{as $t\to\infty$.}
\end{equation}
\par\noindent {\rm (ii)} Under the condition $\alpha=1$ there exists a critical initial velocity $\rv^*$ (that depends only on 
the failure rate), such that if
\begin{equation*}
\|(\bv_1(0),\dots,\bv_k(0))-(\bar{\bv}(0),\dots,\bar{\bv}(0))\|<\rv^*.
\end{equation*}
the statements of {\rm (i)} hold.
\end{theorem}
\begin{remark}The critical velocity $\rv^*$ in {\rm (ii)} above is in fact the expectation of the (random) Fiedler number of 
the
non-colored graph associated to the system. 
We have $\rv^*=0$ when $\lambda=1$ (complete failure), and $\rv^*=k$ when $\lambda=0$ (non-failure), 
see the observations following Remark 2.
\end{remark}
\begin{remark} The case $\lambda=0$ corresponds to the non-failure system, i.e. to the 
Cucker-Smale model introduced in \cite{CS06}.
In case $\alpha=1$ we obtain the condition
$$
\|V(0)\|\leq \rv^*=k\leq \frac1h
$$
according to \eqref{eq:h}. This is similar to the condition in Theorem 3, statement {\rm ii)} in \cite{CS06}. 
(Case $\beta=1/2$ in \cite{CS06} corresponds to case $\alpha=1$ in our framework.)
\end{remark}
A key feature in order to obtain flocking is the connectedness of the graph induced by the agents
of our flock. Within the several ways of quantifying the connectivity of a graph, Cucker and Smale
\cite{CS06} show that a key concept in this situation is the \emph{algebraic connectivity}, also known as the
\emph{Fiedler number} (see \cite{abreu}), that we denote by $\phi$. 
This number is defined as the second smallest eigenvalue of the Laplacian matrix associated to the graph (the first 
eigenvalue always vanishes). The larger the Fiedler number is, the ``more connected'' the graph is.  
A relevant r\^ole will be played also by the Fiedler number of the non-colored graph induced by
the interactions graph, that is defined to have an edge if and only if there is an edge in the original
graph. This second Fiedler number is denoted by $\varphi$. 
In particular, when $\varphi=0$ the graph is not connected, and the complete graph
with $k$ vertices has $\varphi=k$.

Observe that, as the edges of our graph are random, this number also becomes
a random quantity, and the statistical control of its behavior in time gives us the possibility of establishing
our results.

\section{Numerical simulations}
We perform some simulations of the evolution of the system defined by Equations \eqref{eq:cs} 
and \eqref{eq:coefficients}. 
The following graphics correspond to that evolution for some different values of the parameters. 
The graphics show the norm of the system relative velocity, i.e.: the maximum of the $2$-norms of the velocity 
of each agent relative to $v_0$, and the logarithm of that relative velocity.

The initial positions and velocities are sampled from standard normal distribution in $\R^3$, 
and we take $h=1/k$.
\begin{figure}[ht]
\begin{center}
\includegraphics[scale=.9]{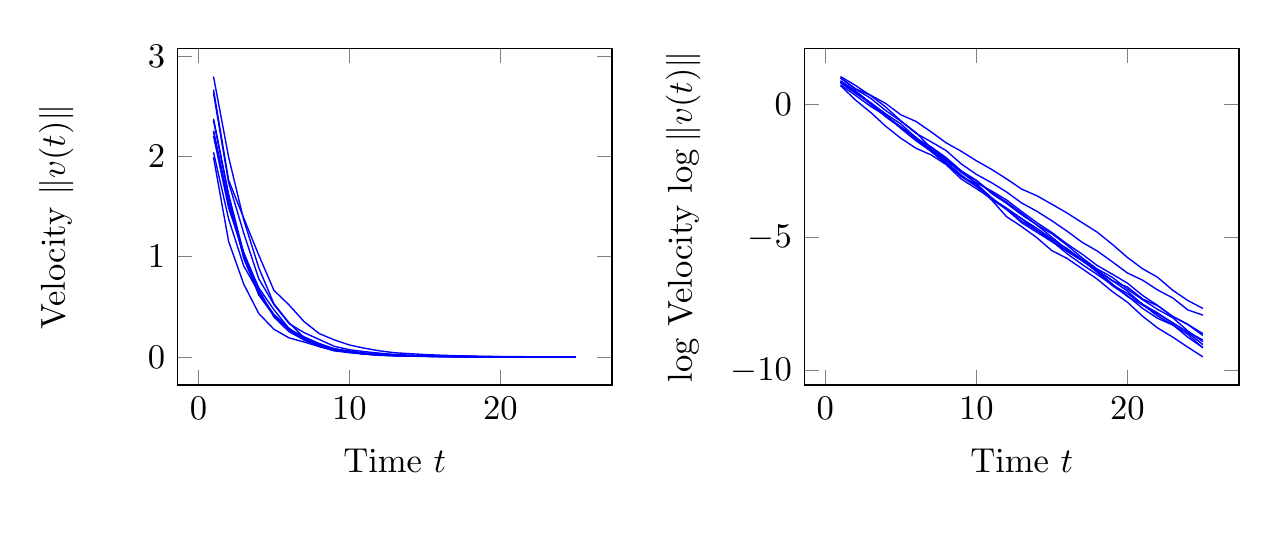}
\caption{Evolution of velocity for $\alpha=0.5$ and $\lambda=0.25$.}
\end{center}
\end{figure}
\begin{figure}[ht]
\begin{center}
\includegraphics[scale=.9]{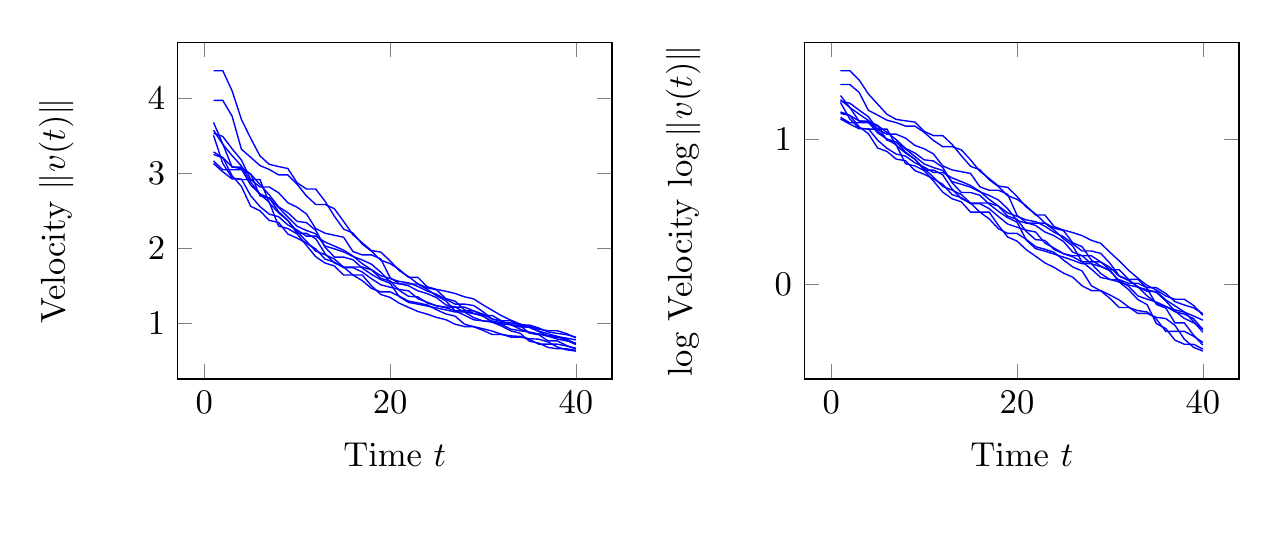}
\caption{Evolution of velocity for $\alpha=0.5$ and $\lambda=0.9$.}
\end{center}
\end{figure}
\begin{figure}[ht]
\begin{center}
\includegraphics[scale=.9]{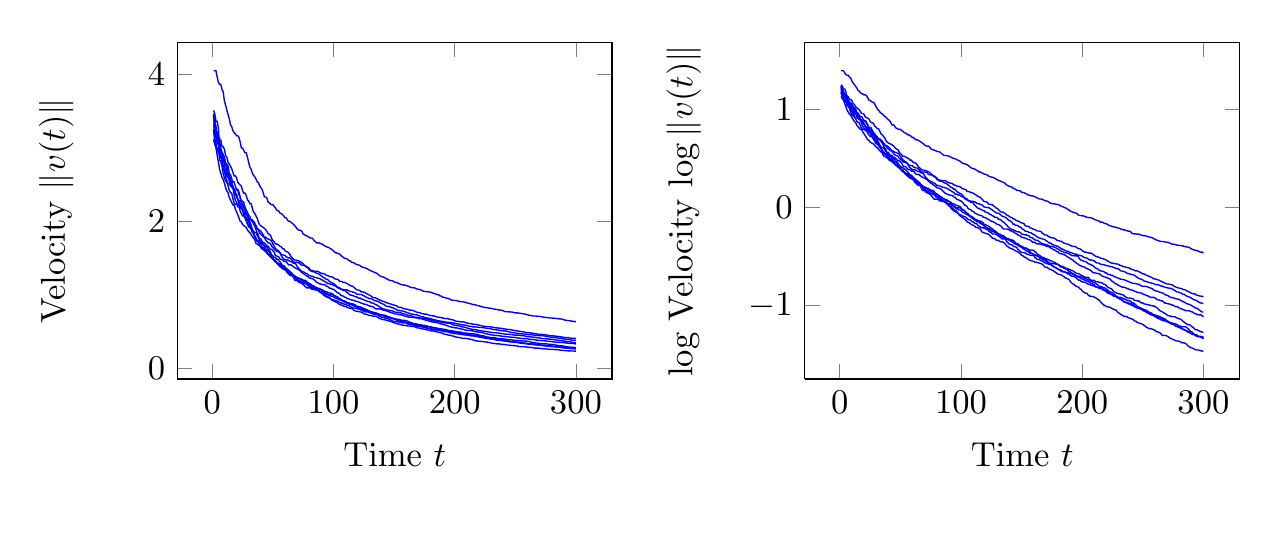}
\caption{Evolution of velocity for $\alpha=1$ and $\lambda=0.9$.}
\end{center}
\end{figure}
It can be observed that the decay seems to be exponential in the sub-linear cases ($\alpha<1$), as in the hierarchical case (see \cite{dm2}) although here our technics to not provide this rate of convergence. The simulation results in the linear case ($\alpha=1$) do not provide clear information about the rate of decay. 

\section{Proof of Theorem \ref{theorem:main}}

%

As was mentioned, we consider the colored graph induced by the system, with vertices $1,2,\dots,k$ corresponding to 
the agents, and edges $a_{ij}(t)$.
We now write  the equations that govern the system in matrix form.
Denote by $\Id$ the $k\times k$ identity matrix, and by $L(t)$ the Laplacian matrix
corresponding to the induced graph  
defined as $L(t)=D(t)-A(t)$ with $A(t)=(a_{ij}(t))$ the incidence matrix of the graph, 
and $D(t)$ the diagonal matrix with $d_{ii}(t)=\sum_ja_{ij}(t)$. 
With this notation, the matrix form of the equations \eqref{eq:cs} are
\begin{align}\label{eq:our}
  \bx(t+h)\;&=\;\bx(t)+h \bv(t),\\
  \bv(t+h)\;&=\;\left(\Id- hL(t)\right)\bv(t),\notag
\end{align}
where the  notation $A\bv$ means that the matrix $A$ is acting on $(\R^3)^k$ by mapping the vector 
$(\bv_1,\dots,\bv_k)$
into the vector $(a_{i1}\bv_1+\cdots+a_{ik}\bv_k)_{1\leq i\leq k}$.

As we have seen in \eqref{eq:vcxc}, the center of mass of the system has constant velocity. 
It is useful then to consider coordinates with respect to this point, introducing the relative position and
velocity of the flock by
\begin{align*}
\rx(t)=(\rx_1(t),\dots,\rx_k(t)) &= (\bx_1(t)-\bar{\bx}(t),\dots,\bx_k(t)-\bar{\bx}(t)),\\
\rv(t)=(\rv_1(t),\dots,\rv_k(t)) &= (\bv_1(t)-\bar{\bv}(0),\dots,\bv_k(t)-\bar{\bv}(0)).
\end{align*}
This change of coordinates is equivalent to the projection on the diagonal space in velocities
introduced by Cucker and Smale in \cite{CS06}. 
A further simplification in the notation is to write $\rx[t]$ and $\rv[t]$ instead of $\rx(th)$ and $\rv(th)$
respectively (and similarly for other time dependent quantities).

With this change of coordinates and notation, the statements of Theorem \ref{theorem:main} are
\begin{equation*}
(\rv_1[t],\dots,\rv_k[t])\to (0,\dots,0),\quad\text{as $t\to\infty$,}
\end{equation*}
and
\begin{equation*}
(\rx_1[t],\dots,\rx_k[t])\to (\hat{\rx}_1,\dots,\hat{\rx}_k),\quad\text{as $t\to\infty$.}
\end{equation*}
It is easy to verify that the relative positions $\rx[t]$ and velocities $\rv[t]$
just introduced follow the same dynamic than the original ones, 
i.e. the system \eqref{eq:our}
holds for $\rx$ and $\rv$ instead of $\bx$ and $\bv$ respectively:
\begin{align*}\label{eq:our}
  \rx[t+1]\;&=\;\rx[t]+h \rv[t],\\
  \rv[t+1]\;&=\;\left(\Id- hL[t]\right)\rv[t].\notag
\end{align*}



We introduce the norm of  a vector $\rx=(\rx_1,\dots,\rx_k)$ in $(\R^3)^k$ by 
$$
\|\rx\|=\sum_{i=1}^k\sum_{\ell=1}^3x_{i\ell}^2.
$$
With this notation the statements \eqref{eq:v-convergence} and \eqref{eq:x-convergence}
reduce to
\begin{equation*}
\|\rv[t]\|\to 0\quad\text{and}\quad\|\rx[t]-\hat{\rx}\|\to 0,\quad\text{as $t\to\infty$}.
\end{equation*}
We further consider the usual operator norm for a matrix $A$ acting as described above by
$$
\|A\|=\sup\{\|Ax\|\colon \|x\|=1\}.
$$
As the Fiedler number of $A$ is the second smallest eigenvalue of $A$, and
the smallest is associated to the eigenvector $(1,\dots,1)$ in $(\R^3)^k$ 
(and our vectors are orthogonal with respect to this diagonal vector), 
from the velocity equation in \eqref{eq:our} we obtain that
\begin{equation}\label{eq:fiedler}
\|\rv[t+1]\|\leq (1-h\phi[t])\|\rv[t]\|,
\end{equation}
where $\phi[t]$ is the Fiedler number of the (random) matrix $\Id- hL[t]$.

First, we observe that the Fiedler number of a colored graph satisfies
$$
\phi\leq \frac{k}{k-1}\min \{d(v): v \in V(G) \}
$$
where $d(v)$ is the degree of the vertex $v$  (see Theorem 2.2 in \cite{mohar}). From this,
as $a_{ij}\leq 1$ for all pairs $i,j$, we obtain $d(v)\leq k-1$ for all $v$, that gives $\phi\leq k$ (the same bound that holds 
for non-colored graphs).
This means, as $h\leq 1/k$ that $0\leq 1-h\phi\leq 1$.

Second, we obtain that the norm of the relative velocity of the flock is decreasing,
giving a linear bound for the norm of the relative position.
\begin{equation*}
\|\rx[t]\|=\|\rx[0]+h(\rv[0]+\cdots+\rv[t-1])\|\leq \|\rx[0]\|+th\|\rv[0]\|.
\end{equation*}
Furthermore, the iteration of \eqref{eq:fiedler} gives
\begin{equation*}
\|\rv[\tau+1]\|\leq \prod_{i=0}^{\tau}(1-h\phi[i])\|\rv[0]\|,
\end{equation*}
that using the equation of the position in \eqref{eq:our}
give us the bound
\begin{eqnarray*}
  \|\rx[\tau]\|  
  &\leq& \|\rx[0]\|+
         h \left(\|\rv[0]\|+\sum_{j=1}^{\tau-1} \|\rv[j]\|\right)\notag\\
  &\leq& \|\rx[0]\|+
        h \|\rv[0]\|\left(1+\sum_{j=1}^{\tau-1} 
          \prod_{i=0}^{j-1}\left(1-h\phi[i]\right)\right).\notag
\end{eqnarray*}
It is crucial then to study the convergence of the series
\begin{equation}\label{eq:s}
S[\tau]=\sum_{j=1}^{\tau-1} 
          \prod_{i=0}^{j-1}\left(1-h\phi[i]\right)
\end{equation}
to obtain an upper bound of the position, that, in its turn, will give us a lower bound
on the Fiedler number.


An important observation is that the connectedness of the coloured graph with incidence matrix
$(a_{ij}[t])$ coincides with the one of the $0$-$1$ graph with incidence matrix $(\zeta^{ij}_t)$
(because both have the same zero and non-zero entries). 
This means that the connectedness of the graph is independent of the position and velocity of the system. 

Let $\varphi[t]$ be the Fiedler number of the non colored graph generated by $(\zeta^{ij}_t)$. 
Then, by Proposition~2 in \cite{CS06}, we have
$$
\phi[t]\geq\varphi[t]\mu[t]
$$
with 
$$\mu[t]  = \min\{a_{ij}[t]\colon a_{ij}[t] >0, i\neq j\}.$$

Note that when the graph is not connected we have $\phi[t]=\varphi[t]=0$.


We now prove that  $\mu[t] \geq \frac{A}{B + t^\alpha}$ for some constants $A$ and $B$. 
For this, we rely on the inequality $(a+b)^\alpha\leq a^\alpha+b^\alpha$ that holds for $a\geq0$, $b\geq 0$ and $0\leq
\alpha\leq 1$.
Applying this inequality with $a=1+\|x[0]\|$ and $b=h\|v[0]\|t$ we obtain
$$
\mu[t]\geq\frac{1/(h\|v[0]\|)^\alpha}{\left(\frac{1+\|x[0]\|}{h\|v[0]\|}\right)^\alpha+t^\alpha}=:\frac{A}{B+t^\alpha}
$$
for all $t$, that gives 


$$
\phi[t]\geq\varphi[t]\frac{A}{B+t^\alpha}=:\nu[t].
$$

With this result we obtain the following bound for our sum in \eqref{eq:s}.
Denote $\bar{\varphi}=\mean{\varphi[t]}$ and note that $\bar{\varphi}>0$ since $\varphi$ is a nonnegative and 
not identically zero random variable. 
As the quantities $\{\nu[t]\}$ form a sequence of independent (non identically distributed) 
random variables, and
$$
\mean(1-h\nu[t])=1-h\bar{\varphi}{\frac{A}{B+t^{\a}}}
$$
Furthermore
\begin{align*}
\mean S[\tau]&\leq \sum_{j=1}^{\tau-1}\prod_{i=1}^{j}\mean\left(1-h\nu[i]\right)=
\sum_{j=1}^{\tau-1}\prod_{i=1}^{j}
\left(1-{\bar{\varphi} h \frac{A}{B+i^{\a}}}\right).\\
\end{align*}
Assume $\alpha<1$. 
For each summand above, as $\log(1-x)\leq -x$,  we have
\begin{align}
\prod_{i=1}^{j}
\left(1-{\bar{\varphi} h \frac{A}{B+i^{\a}}}\right)&
=\exp\left(\sum_{i=1}^j\log
\left(1-{\bar{\varphi} h \frac{A}{B+i^{\a}}}\right)\right)
\leq 
\exp\left(-
\sum_{i=1}^j{\bar{\varphi} h \frac{A}{B+i^{\a}}}\right)\notag\\
&\leq
\exp\left(-\bar{\varphi} h A
\sum_{i=1}^j{\frac{1}{i^{\a}}}\right)\leq 
\exp\left(-\gamma j^{1-\a}\right),\label{eq:term}
\end{align}
where $\gamma=(\bar{\varphi} h A)/(1-\alpha)$.
Finally we observe that a series with general term given by \eqref{eq:term} is convergent,
and this implies the convergence of $\mean(S[\tau])$ as $\tau\to\infty$ to a finite limit.
As the series $S[\tau]$ itself is increasing with $\tau$, we obtain that there exists
an almost sure finite limit
\begin{equation}\label{eq:limit}
S=\lim_{\tau\to\infty}S[\tau]<\infty.
\end{equation}

The case $\alpha=1$ is treated separately. 
We have to modify the final bound in \eqref{eq:term}, in this case we have
\begin{align*}
\prod_{i=1}^{j}
\left(1-{\bar{\varphi} h \frac{A}{B+i}}\right)&
=\exp\left(\sum_{i=1}^j\log
\left(1-{\bar{\varphi} h \frac{A}{B+i}}\right)\right)
\leq 
\exp\left(-\bar{\varphi} h A
\sum_{i=1}^j{ \frac{1}{B+i}}\right)\notag\\
&\leq\exp\left(-\bar{\varphi}hA\log\left(\frac{B+j+1}{B+1}\right)\right)
=\left(\frac{B+1}{B+j+1}\right)^{\bar{\varphi}hA}.
\end{align*}

Now, a series with this general term converges when $\bar{\varphi}hA>1$, that is $\|v[0]\|<\bar{\varphi}$.
Under this condition we obtain \eqref{eq:limit}.


From the convergence obtained in \eqref{eq:limit} (in both cases $\alpha<1$ and $\alpha=1$)
for the series defined in \eqref{eq:s}, we get that the series of the velocities norm is convergent, i.e.
$$
\sum_{j=0}^{\tau-1} \|v[j]\|\leq \|v[0]\|(1+S)
$$
This implies the convergence stated in \eqref{eq:v-convergence}. It also implies that the series of the velocities 
converges itself, i.e. there exists 
$$
\hat{\rv}=\sum_{j=0}^{\infty} \rv[j].
$$
This fact gives the convergence for the positions of the system:
\begin{eqnarray*}
  \rx[\tau]  
  &=&\rx[0]+\sum_{j=0}^{\tau-1}
     \left(\rx[j+1]-\rx[j]\right)
  = \rx[0]+
         h \sum_{j=0}^{\tau-1} \rv[j]\notag\\
  &\to& \rx[0]+
         h \hat{\rv}=:\hat{\rx},
\end{eqnarray*}
concluding the proof of the Theorem.

\end{document}